\documentclass[12pt]{article}

\usepackage{fullpage}
\usepackage{amssymb}
\usepackage{amsmath}
\usepackage[usenames]{color}
\usepackage{graphicx}
\usepackage[colorlinks=true,
linkcolor=webgreen,
filecolor=webbrown,
citecolor=webgreen]{hyperref}

\definecolor{webgreen}{rgb}{0,.5,0}
\definecolor{webbrown}{rgb}{.6,0,0}
\definecolor{red}{rgb}{1,0,0}

\usepackage{color}
\newtheorem{theorem}{Theorem}

\newtheorem{example}[theorem]{Example}

\newtheorem{remark}[theorem]{Remark}

\begin{document}
	
	\title{Partitions with unique largest part and their generating functions}
	\date{}
	\author{H. Kaur$^1$ and M. Rana$^2$\\
	$^1\,$Department of Mathematics, Chandigarh University,\\ Mohali-140413, Punjab,\\
	$^2\,$Department of Mathematics, Thapar Institute of Engineering and Technology,\\ Patiala-147004, Punjab, India\\
\url{kaur.harman196@gmail.com}\\
\url{mrana@thapar.edu}}
	\maketitle
	
	\begin{center}
			{Abstract}
			\end{center}
			The paper introduce a new type of partitions where the largest part appears exactly once, and the remaining parts constitute a partition of that largest part. We derive the generating function associated with these partitions and subsequently explore several variations, providing the corresponding generating functions for each variant. 

\section{Introduction}
	Partitions have been the object of interest to mathematicians for many years. A \textit{Partition} of a positive integer $n$ is a way of expressing it as a sum of positive integer parts, generally the parts are listed in descending order, from the largest ($\ell$) to the smallest ($s$). The generating function to count the number of partitions $p(n)$ of $n$ is $$\sum_{n \geq 0} p(n)q^n = \frac{1}{(q;q)_{\infty}}$$
	where 
	$$
	(t; q)_n = \begin{cases} 
		1 &  \text{ for }n=0, \\
		(1-t)(1-tq)(1-tq^2)(1-tq^3)\cdots (1-tq^{n-1}) & \text{ for } n > 0 
	\end{cases}
	$$
	is the $q-$shifted factorial and
	$$	(t; q)_{\infty} = \lim_{n\to\infty} (t;q)_n, \text {~~where~~} |q|<1.$$
	
	For a detailed understanding of integer partitions, one can refer to \cite{Aggarwal, Andrews}. The study of restricted partitions dates back to Euler, who established a remarkable result: the number of partitions of any integer into odd parts, $o(n)$, is equal to the number of partitions into distinct parts or strict partitions, $d(n)$, the corresponding sequence of partitions is listed in On Line Encyclopedia of Integer Sequence (OEIS) with number A000009. The OEIS is a vast database of well-defined sequences of integers, many of which originate from partitions theory. Analogous to the results of Euler, many appealing results on partitions can be found in \cite{Gleissberg, Hirschhorn, Merca1, Schur}.\\\\
	In this paper, we discuss the partition, denoted by $\Lambda$ where the largest part appears exactly once, and the remaining parts constitute a partition of that largest part. For example, let $n=10$, then some relevant partitions, $\Lambda$, satisfying the condition are: 
	$$5+4+1,5+3+2, 5+2+2+1$$
	However, several other partitions of $n=10,$ such as $9+1, 8+2, 7+3,$ do not satisfy the condition. Since the sum of all parts other than the largest part always equals the largest part ($\ell$), the total sum is always an even number ($2\ell$).  We provide the generating function corresponding to these partitions and to some associated restricted partitions analogously. The integer sequences corresponding to partitions described in Theorems \ref{theorem1}-\ref{theorem3a} are listed in OEIS. However, the partition sequences obtained in Theorems \ref{theorem4}-\ref{theorem5} are new.
		\begin{theorem} \label{theorem1}
		The generating function to count the number of partitions $\Lambda$ is 
		$$	\sum_{n=0}^{\infty}\rho(n)q^n = \frac{1}{(q^2;q^2)_{\infty}}- \frac{1}{(1-q^2)}.$$
	\end{theorem}
	
	We now consider some variations of the partition $\Lambda$, motivated by Euler’s theorem on partitions of $n$ into distinct or odd parts.  The remaining parts, other than the largest part, comes with some restrictions. Analogously, we obtained their generating functions. We use the notations as follows:
	\begin{align*}
		\Lambda_d: & \text{ partition } \Lambda \text{ where the parts are distinct,} \\
		\Lambda_o: & \text{ partition } \Lambda \text{ where the parts are odd, } \\
		\Lambda_{od}: & \text{ partition } \Lambda \text{ where the parts are odd and distinct. }
	\end{align*}
	\begin{theorem}  \label{theorem2}
		The generating function to count the number of partitions $\Lambda_d$ is
		$$\sum_{n=0}^{\infty}\rho_d(n)q^n = {(-q^2;q^2)_{\infty}}- \frac{1}{(1-q^2)}.$$
	\end{theorem}
	\begin{theorem}  \label{theorem3}
		The generating function to count the number of partitions $\Lambda_o$ is
		$$\sum_{n=0}^{\infty}\rho_o(n)q^n = \frac{1}{(q^2;q^4)_{\infty}}- \frac{q^2}{(1-q^4)}-1.$$
	\end{theorem}
	\begin{theorem}  \label{theorem3a}
		The generating function to count the number of partitions $\Lambda_{od}$ is
		$$\sum_{n=0}^{\infty}\rho_{od}(n)q^n = {(-q^2;q^4)_{\infty}}- \frac{q^2}{(1-q^4)}-1.$$
	\end{theorem}
	
	Further inspired by Schur theorem which states that the number of partitions into distinct parts  $\equiv \pm 1 \pmod3$ is same as the number of partitions into parts  $\equiv \pm 1 \pmod6$. We consider the partitions $\Lambda$ with distinct partss and $\equiv \pm 1 \pmod3$ except the largest part, denoted by $\Lambda_3$. 
	\begin{theorem}  \label{theorem4}
		The generating function to count the number of partitions $\Lambda_3$ is
		$$\sum_{n=0}^{\infty}\rho_3(n)q^n = (-q^2,-q^4;q^6)_{\infty} - \frac{q^2}{1-q^2}+ \frac{q^6}{1-q^6}.$$
	\end{theorem}
	Analogously, the generating function for the partitions $\Lambda$ with partss  $\equiv \pm 1 \pmod6$ except the largest part, denoted by $\Lambda_6$ is given in the next theorem. 
	\begin{theorem}  \label{theorem5}
		The generating function to count the number of partitions $\Lambda_6$ is
		$$\sum_{n=0}^{\infty}\rho_6(n)q^n =\frac{1}{(q^2,q^{10};q^{12})_{\infty}} - \frac{(q^2+q^{10})}{1-q^{12}}.$$
	\end{theorem}
	\section{Proof of the theorems}
\textbf{Proof of Theorem \ref{theorem1}}
		The partitions $\Lambda,$ with the largest part and rest of the parts forms the partitions of largest part, can be generated as: 
		\begin{align*}
			\sum_{n=0}^{\infty}\rho(n)q^n &= q^2 \cdot q^{1+1} +q^3 \cdot (q^{2+1}+q^{1+1+1} ) + q^4 \cdot (q^{3+1}+q^{2+2}+q^{2+1+1}+q^{1+1+1+1})\\
			& \hspace{2mm}+q^5 \cdot (q^{4+1}+q^{3+2}+q^{3+1+1}+q^{2+2+1}+q^{2+1+1+1}+q^{1+1+1+1+1})+ \cdots \\
			\sum_{n=0}^{\infty}\rho(n)q^n &= \sum_{n=0}^{\infty}q^n (p(n)-1)q^n\\
			&= \sum_{n=0}^{\infty} (p(n)-1)q^{2n}\\
			&= \sum_{n=0}^{\infty} p(n)q^{2n} - \sum_{n=0}^{\infty} q^{2n}\\
			&= \frac{1}{(q^2;q^2)_{\infty}}- \frac{1}{(1-q^2)}.
		\end{align*}
\textbf{Proof of Theorem \ref{theorem2}}
		Consider partitions $\Lambda_d,$ with the largest part and rest of the parts are distinct that represents the partition of the corresponding largest part, it can be generated as:
		\begin{align*}
			\sum_{n=0}^{\infty}\rho_d(n)q^n &=q^3 \cdot (q^{2+1} ) + q^4 \cdot (q^{3+1})+q^5 \cdot (q^{4+1}+q^{3+2})+q^6 \cdot (q^{5+1}+q^{4+2}+q^{3+2+1})\\
			& \hspace{2mm}+ q^7 \cdot (q^{6+1}+q^{5+2}+q^{4+3}+q^{4+2+1})+\cdots .
		\end{align*}
		Let $d(n)$ represents the number of distinct partitions of $n$, then
		\begin{align*}
			\sum_{n=0}^{\infty}\rho_d(n)q^n &= \sum_{n \geq 3}q^n (d(n)-1)q^n\\
			&= \sum_{n \geq 3}^{\infty} (d(n)-1)q^{2n}\\
			&= \sum_{n\geq 3} d(n)q^{2n} - \sum_{n \geq 3} q^{2n}\\
			 &= \sum_{n\geq 0} d(n)q^{2n} - \sum_{n \geq 0} q^{2n}\\
			\sum_{n=0}^{\infty}\rho_d(n)q^n	
			&= {(-q^2;q^2)_{\infty}}- \frac{1}{(1-q^2)}.
		\end{align*}
\textbf{Proof of Theorem \ref{theorem3}}
		The generating function for partitions with the largest part and rest of the parts are odd that constitutes the partition of the corresponding largest part leads to
		\begin{align*}
			\sum_{n=0}^{\infty}\rho_o(n)q^n &=q^2 \cdot (q^{1+1})+q^3 \cdot (q^{1+1+1} ) + q^4 \cdot (q^{3+1}+q^{1+1+1})+\\
			&\hspace{0.7cm}q^5 \cdot (q^{3+1+1}+q^{1+1+1+1+1})+
			q^6 \cdot (q^{5+1}+q^{3+3}+q^{3+1+1+1})+\\
			& \hspace{1cm} q^7 \cdot (q^{5+1+1}+q^{3+3+1}+q^{3+1+1+1+1}+q^{1+1+1+1+1+1+1})+\cdots
		\end{align*}
		Let $o(n)$ represents the number of odd partitions of $n$, then
		\begin{align*}
			\sum_{n=0}^{\infty}\rho_o(n)q^n &= \sum_{\substack{{n \geq 2}\\{n-even}}} o(n)q^{2n} + \sum_{\substack{{n \geq 3}\\{n-odd}}} (o(n)-1)q^{2n}\\
			&= \sum_{n\geq 2} o(n)q^{2n} - \sum_{\substack{{n \geq 3}\\{n-odd}}} q^{2n}\\
		&= \sum_{n\geq 0} o(n)q^{2n} - \sum_{\substack{{n \geq 1}\\{n-odd}}} q^{2n}-(o(0)q^0 +o(1)q^2)+q^2\\
			\end{align*}
		\begin{align*}
		\sum_{n=0}^{\infty}\rho_o(n)q^n	&= \sum_{n\geq 0} o(n)q^{2n} - \sum_{\substack{{n \geq 1}\\{n-odd}}} q^{2n}-(o(0)q^0 +o(1)q^2)+q^2\\
			&= \frac{1}{(q^2;q^4)_{\infty}}- \frac{q^2}{(1-q^4)}-1
		\end{align*}
\textbf{Proof of Theorem \ref{theorem3a}}
		The generating function for partitions $\Lambda$ with the parts (apart from the largest part) are distinct and odd that constitutes the largest part leads to
		\begin{align*}
			\sum_{n=0}^{\infty}\rho_{od}(n)q^n &=q^4 \cdot q^{3+1}+q^6 \cdot q^{5+1} +q^8 \cdot (q^{7+1}+q^{5+3}) +q^9 \cdot (q^{5+3+1})+\\
			&\hspace{1cm}q^{10} \cdot (q^{9+1}+q^{7+3})+q^{11} \cdot (q^{7+3+1})+ q^{12} \cdot (q^{11+1}+q^{9+3}+q^{7+5})+\\
			& \hspace{1cm}q^{13} \cdot (q^{9+3+1}+q^{7+5+1})+\cdots
		\end{align*}
		Let $od(n)$ represents the number of odd and distinct partitions of $n$, then
		\begin{align*} 
			\sum_{n=0}^{\infty}\rho_{od}(n)q^n&= \sum_{\substack{{n \geq 1}\\{n \text{ even}}}}q^n od(n)q^n + \sum_{\substack{{n \geq 1}\\{n \text{ odd}}}}q^n (od(n)-1)q^n\\
			&= \sum_{n \geq 1} od(n)q^{2n} - \sum_{\substack{{n \geq 1}\\{n \text{ odd}}}} q^{2n}\\
			&= \sum_{n \geq 0} od(n)q^{2n} - \frac{q^2}{1-q^4}-1
		\end{align*}
		\begin{align*}
			\sum_{n=0}^{\infty}\rho_{od}(n)q^n	
			&= (-q^2;q^4)_{\infty} - \frac{q^2}{1-q^4}-1.
		\end{align*}
	\begin{example} In this example, we present the relevant partitions for $n=8$ corresponding to $\Lambda$ and its variations, as discussed above,  along with the associated sequence numbers and combinatorial interpretations as mentioned in the OEIS :\\\\
	\footnotesize{
		\begin{table}[h] \label{table1}\begin{center}
				
				\bigskip
				{\renewcommand{\arraystretch}{1.2}
					\begin{tabular}{|c|c|c|c|}
						\hline
						For $n=8$ & Relevant partitions  & OEIS sequence& OEIS interpretation \\
						{} & & {no.} & {}\\
						\hline
						$\rho(8)=4$ &  $4+3+1, 4+2+2, $& A000065 & $-1$ + \# partitions of $n$\\
						{} & $ 4+2+1+1,4+1+1+1+1$ & {} & {}\\
						\hline
						$\rho_d(8)=1$ &  $4+3+1$& A111133  & \# partitions of $n$ into \\
						{} & {} & {} & at least two distinct parts\\
						\hline
						$\rho_o(8)=2$ &  $4+3+1,4+1+1+1+1$& 	A357456 & \# partitions of $n$ into\\
						{} & {} & {} &  two or more odd parts\\ 
						\hline
						$\rho_{od}(10)=2$ & $4+3+1$& 	A357457  & \# partitions of $n$ into\\
						{} & {} & {} &  two or more distinct odd parts\\
						\hline
					\end{tabular}
			}\end{center}
		\end{table}
	}
	
	\end{example}

\noindent \textbf{Proof of Theorem \ref{theorem4}}
		The generating function for partitions $\Lambda$ with the parts (apart from the largest part) are distinct and $ \equiv\pm 1 \pmod3$ that constitutes the largest part leads to
		\begin{align*}
			\sum_{n=0}^{\infty}\rho_3(n)q^n &=q^3 \cdot q^{1+2}+q^5 \cdot q^{4+1} +q^6 \cdot (q^{5+1}+q^{4+2}) +q^7 \cdot (q^{5+2}+q^{4+2+1})+\\
			&\hspace{1cm}q^8 \cdot (q^{7+1}+q^{5+2+1})+q^9 \cdot (q^{8+1}+q^{7+2}+q^{5+4})+ \\
			& \hspace{1cm}q^{10} \cdot (q^{8+2}+q^{7+2+1}+q^{5+4+1})+q^{11} \cdot (q^{10+1}+q^{8+2+1}+q^{7+1}+q^{5+4+2})\cdots\\
		\end{align*}
		Let $dp_3(n)$ represents the number of distinct partitions of $n$ with parts $ \equiv\pm 1 \pmod3$, then
		\begin{align*} 
			\sum_{n=0}^{\infty}\rho_3(n)q^n&= \sum_{\substack{{n \geq 3}\\{3 | n}}}q^n dp_3(n)q^n + \sum_{\substack{{n \geq 1}\\{3 \nmid n}}}q^n (dp_3(n)-1)q^n\\
			&= \sum_{n \geq 1} dp_3(n)q^{2n} - \sum_{\substack{{n \geq 1}\\{3 \nmid n}}} q^{2n}\\
			&= \sum_{n \geq 1} dp_3(n)q^{2n} - \frac{q^2}{1-q^2}+ \frac{q^6}{1-q^6}\\
			&= (-q^2,-q^4;q^6)_{\infty} - \frac{q^2}{1-q^2}+ \frac{q^6}{1-q^6}.
		\end{align*}

\textbf{Proof of Theorem \ref{theorem5}}
		The generating function for partitions $\Lambda$ with the partss (apart from the largest part) are $ \equiv\pm 1 \pmod6$ that constitutes the largest part leads to
		\begin{align*}
			\sum_{n=0}^{\infty}\rho_6(n)q^n &=q^2 \cdot q^{1+1}+q^3 \cdot q^{1+1+1} +q^4 \cdot (q^{1+1+1+1}) +q^5 \cdot (q^{1+1+1+1+1})\\
			&+q^6 \cdot (q^{5+1}+q^{1+1+1+1+1+1})+q^7 \cdot (q^{5+1+1}+q^{1+1+1+1+1+1+1})\\
			&+ q^{8} \cdot (q^{7+1}+q^{5+1+1+1}+q^{1+1+1+1+1+1+1+1})\\
			&+q^{9} \cdot (q^{7+1+1}+q^{5+1+1+1+1}+q^{1+1+1+1+1+1+1+1+1})\cdots
		\end{align*}
		Let $p_6(n)$ represents the number of partitions of $n$ with parts $ \equiv\pm 1 \pmod6$, then
		\begin{align*}
			\sum_{n=0}^{\infty}\rho_6(n)q^n &= \sum_{\substack{{n \geq 2}\\{n \not\equiv \pm 6}}}q^n p_6(n)q^n + \sum_{\substack{{n \geq 1}\\{n \equiv \pm 6}}}q^n (p_6(n)-1)q^n
		\end{align*}
		
		\begin{align*}
			&= \sum_{n \geq 1} p_6(n)q^{2n} - \sum_{\substack{{n \geq 1}\\{n \equiv \pm 6}}} q^{2n}\\
			&= \frac{1}{(q^2,q^{10};q^{12})_{\infty}} - \frac{(q^2+q^{10})}{1-q^{12}}.
		\end{align*}

\begin{remark}

The partition sequences obtained in Theorem 5-6 are new, given as: \\
\noindent Sequence 1: $<1,1,1,2,2,2,3,3,4,6,6,7,9,9,11,14,15,17,20,22,25,30,33,37,42,\ldots>$\\
\noindent Sequence 2: $<1,1,1,1,1,2,2,3,3,4,4,6,6,8,9,10,11,14,15,18,20,23,25,30,33,\ldots>$\\
\noindent Here, the former sequence 1 represents the coefficients $\rho_3(n)$ due to Theorem \ref{theorem4} while the later sequence represents the coefficients $\rho_6(n)$ due to Theorem \ref{theorem5}.

\end{remark}

\end{document}